\documentclass[11pt]{article}
\usepackage[]{amsmath,amssymb}

\newtheorem{theorem}{Theorem}[section]
\newtheorem{lemma}[theorem]{Lemma}
\newtheorem{proposition}[theorem]{Proposition}
\newtheorem{corollary}[theorem]{Corollary}
\newtheorem{exAux}[theorem]{Example}
\newenvironment{example}{\begin{exAux} \rm}{\end{exAux}}
\newtheorem{Def}[theorem]{Definition}
\newenvironment{definition}{\begin{Def} \rm}{\end{Def}}
\newtheorem{Note}[theorem]{Note}
\newenvironment{note}{\begin{Note} \rm}{\end{Note}}
\newtheorem{Rem}[theorem]{Remark}

\newtheorem{Ass}[theorem]{Assumption}

\newenvironment{proof}{\medskip\noindent{\bf Proof.\ }}{\qed\medskip}

\newcommand{\qed}{\hfill\mbox{$\Box$\qquad\qquad}}

\newcommand{\Mat}[1]{\text{\rm Mat}_{#1}(\mathbb{K})}

\begin{document}
\thispagestyle{empty}

\begin{center}
\LARGE \bf
\noindent
Matrix units associated with the split basis \\
of a Leonard pair
\end{center}

\smallskip

\begin{center}
\Large
Kazumasa Nomura and Paul Terwilliger
\end{center}

\smallskip

\begin{quote}
\small 
\begin{center}
\bf Abstract
\end{center}
Let $\mathbb{K}$ denote a field, and let $V$ denote a vector
space over $\mathbb{K}$ with finite positive dimension.
We consider a pair of linear transformations
$A:V \to V$ and $A^*:V \to V$ that satisfy (i), (ii) below:
\begin{itemize}
\item[(i)] There exists a basis for $V$ with respect to which the
matrix representing $A$ is irreducible tridiagonal and the matrix
representing $A^*$ is diagonal.
\item[(ii)] There exists a basis for $V$ with respect to which the
matrix representing $A^*$ is irreducible tridiagonal and the matrix
representing $A$ is diagonal.
\end{itemize}
We call such a pair a {\em Leonard pair} on $V$.
It is known that there exists a basis for $V$ with respect to which
the matrix representing $A$ is lower bidiagonal and the matrix
representing $A^*$ is upper bidiagonal. 
In this paper we give some formulae involving the
matrix units associated with this basis. 
\end{quote}

\section{Leonard pairs and Leonard systems}

We begin by recalling the notion of a Leonard pair.
We will use the following terms.
A square matrix $X$ is said to be {\em tridiagonal}
whenever each nonzero entry lies on either the diagonal, the subdiagonal,
or the superdiagonal. Assume $X$ is tridiagonal.
Then $X$ is said to be {\em irreducible}
whenever each entry on the subdiagonal is nonzero and each entry on
the superdiagonal is nonzero.
We now define a Leonard pair.
For the rest of this paper $\mathbb{K}$ will denote a field.

\medskip

\begin{definition}  \cite{T:Leonard}  \label{def:LP}
Let $V$ denote a vector space over $\mathbb{K}$ with finite positive
dimension.
By a {\em Leonard pair} on $V$ we mean an ordered pair $(A,A^*)$,
where $A:V \to V$ and $A^*:V \to V$ are linear transformations
that satisfy (i), (ii) below:
\begin{itemize}
\item[(i)] There exists a basis for $V$ with respect to which the
matrix representing $A$ is irreducible tridiagonal and the matrix
representing $A^*$ is diagonal.
\item[(ii)] There exists a basis for $V$ with respect to which the
matrix representing $A^*$ is irreducible tridiagonal and the matrix
representing $A$ is diagonal.
\end{itemize}
\end{definition}

\begin{note}
It is a common notational convention to use $A^*$ to represent the
conjugate-transpose of $A$. We are {\em not} using this convention.
In a Leonard pair $(A,A^*)$ the linear transformations $A$ and
$A^*$ are arbitrary subject to (i) and (ii) above.
\end{note}

\medskip

We refer the reader to 
\cite{H},
\cite{N:aw}, 
\cite{NT:balanced}, \cite{NT:formula}, \cite{NT:det},
\cite{P}, \cite{T:sub1}, \cite{T:sub3}, \cite{T:Leonard},
\cite{T:24points}, \cite{T:conform}, \cite{T:intro},
\cite{T:intro2}, \cite{T:split}, \cite{T:array}, \cite{T:qRacah},
\cite{T:survey}, \cite{TV}, \cite{V}
for background on Leonard pairs.
We especially recommend the survey \cite{T:survey}.
See \cite{AC}, \cite{AC2}, \cite{ITT}, \cite{IT:shape},
\cite{IT:uqsl2hat}, \cite{IT:non-nilpotent}, \cite{ITW:equitable}, 
\cite{N:refine}, \cite{N:height1},
\cite{T:qSerre}, \cite{T:Kac-Moody} for related topics.

\medskip

It is known that there exists a basis for $V$ with respect to which
the matrix representing $A$ is lower bidiagonal
and the matrix representing $A^*$ is upper bidiagonal 
\cite[Theorem 3.2]{T:Leonard}.
In this paper we give some formulae involving
the matrix units associated with this basis.
We also display some related formulae involving the primitive
idempotents of $A$ and $A^*$,
which might be of independent interest.

\medskip

When working with a Leonard pair, it is convenient to consider a closely
related object called a {\em Leonard system}. 
To prepare for our definition
of a Leonard system, we recall a few concepts from linear algebra.
Let $d$ denote a nonnegative integer and let
$\Mat{d+1}$ denote the $\mathbb{K}$-algebra consisting of all $d+1$ by
$d+1$ matrices that have entries in $\mathbb{K}$. 
We index the rows and  columns by $0, 1, \ldots, d$. 
For the rest of this paper we let $\cal A$ denote a $\mathbb{K}$-algebra 
isomorphic to $\Mat{d+1}$. 
Let $V$ denote a simple $\cal A$-module. We remark that $V$ is unique
up to isomorphism of $\cal A$-modules, and that $V$ has dimension $d+1$.
Let $v_0, v_1, \ldots, v_d$ denote a basis for $V$.
For $X \in {\cal A}$ and $Y \in \Mat{d+1}$, we say 
{\em $Y$ represents $X$ with respect to} $v_0, v_1, \ldots, v_d$
whenever $X v_j = \sum_{i=0}^d Y_{ij}v_i$ for $0 \leq j \leq d$.
For $A \in \cal A$ we say $A$ is {\em multiplicity-free}
whenever it has $d+1$ mutually distinct eigenvalues in $\mathbb{K}$. 
Assume $A$ is multiplicity-free. 
Let $\theta_0, \theta_1, \ldots, \theta_d$ denote an ordering 
of the eigenvalues of $A$, and for $0 \leq i \leq d$ put
\begin{equation}        \label{eq:defEi}
    E_i = \prod_{\stackrel{0 \leq j \leq d}{j\neq i}}
             \frac{A-\theta_j I}{\theta_i - \theta_j},
\end{equation}
where $I$ denotes the identity of $\cal A$. 
 We observe
(i) $AE_i = \theta_i E_i$ $(0 \leq i \leq d)$;
(ii) $E_i E_j = \delta_{i,j} E_i$ $(0 \leq i,j \leq d)$;
(iii) $\sum_{i=0}^{d} E_i = I$;
(iv) $A = \sum_{i=0}^{d} \theta_i E_i$.
We call $E_i$ the {\em primitive idempotent} of $A$ associated with
$\theta_i$. 

By a {\em Leonard pair in $\cal A$} we mean an ordered pair of elements
in $\cal A$ that act on $V$ as a Leonard pair in the sense of
Definition \ref{def:LP}.
We now define a Leonard system.

\medskip

\begin{definition}  \cite{T:Leonard}     \label{def:LS}
By a {\em Leonard system} in $\cal A$ we mean a sequence
\[
   (A; \{E_i\}_{i=0}^d; A^*; \{E^*_i\}_{i=0}^d)
\]
that satisfies (i)--(v) below.
\begin{itemize}
\item[(i)] Each of $A$, $A^*$ is a multiplicity-free element in $\cal A$.
\item[(ii)] $E_0, E_1, \ldots, E_d$ is an ordering of the
   primitive idempotents of $A$.
\item[(iii)] $E^*_0, E^*_1, \ldots, E^*_d$ is an ordering of the
   primitive idempotents of $A^*$.
\item[(iv)] For $0 \leq i,j \leq d$, 
\[
   E_i A^* E_j =
    \begin{cases}  
        0 & \text{\rm if $|i-j|>1$},  \\
        \neq 0 & \text{\rm if $|i-j|=1$}.
    \end{cases}
\]
\item[(v)] For $0 \leq i,j \leq d$, 
\[
   E^*_i A E^*_j =
    \begin{cases}  
        0 & \text{\rm if $|i-j|>1$},  \\
        \neq 0 & \text{\rm if $|i-j|=1$}.
    \end{cases}
\]
\end{itemize}
\end{definition}

\medskip

Leonard systems are related to Leonard pairs as follows.
Let $(A; \{E_i\}_{i=0}^d;A^*; \{E^*_i\}_{i=0}^d)$ denote a Leonard system
in $\cal A$. Then $(A,A^*)$ is a Leonard pair in $\cal A$
\cite[Section 3]{T:qRacah}.
Conversely, suppose $(A,A^*)$ is a Leonard pair in $\cal A$.
Then each of $A,A^*$ is multiplicity-free \cite[Lemma 1.3]{T:Leonard}.
Moreover there exists an ordering $E_0,E_1,\ldots,E_d$ of the
primitive idempotents of $A$, and 
there exists an ordering $E^*_0,E^*_1,\ldots,E^*_d$ of the
primitive idempotents of $A^*$, such that
$(A; \{E_i\}_{i=0}^d; A^*; \{E^*_i\}_{i=0}^d)$
is a Leonard system in $\cal A$ \cite[Lemma 3.3]{T:qRacah}.

\section{The $D_4$ action}

For the rest of this paper, let 
\begin{equation}
\Phi=(A; \{E_i\}_{i=0}^d; A^*; \{E^*_i\}_{i=0}^d)
\end{equation}
denote a Leonard system in $\cal A$.

\medskip

Observe that each of the following is a Leonard system in $\cal A$:
\begin{eqnarray*}
\Phi^{*}  &:=& 
       (A^*; \{E^*_i\}_{i=0}^d; A; \{E_i\}_{i=0}^d),  \\
\Phi^{\downarrow} &:=&
       (A; \{E_i\}_{i=0}^d; A^*; \{E^*_{d-i}\}_{i=0}^d), \\
\Phi^{\Downarrow} &:=&
       (A; \{E_{d-i}\}_{i=0}^d; A^*; \{E^*_{i}\}_{i=0}^d).
\end{eqnarray*}
Viewing $*$, $\downarrow$, $\Downarrow$ as permutations on the set of
all Leonard systems in $\cal A$,
\begin{equation}    \label{eq:relation1}
*^2 = \downarrow^2 = \Downarrow^2 = 1,
\end{equation}
\begin{equation}    \label{eq:relation2}
\Downarrow * = * \downarrow, \quad
\downarrow * = * \Downarrow, \quad
\downarrow \Downarrow = \Downarrow \downarrow.
\end{equation}
The group generated by symbols $*$, $\downarrow$, $\Downarrow$ subject
to the relations (\ref{eq:relation1}), (\ref{eq:relation2}) is the
dihedral group $D_4$. We recall $D_4$ is the group of symmetries of a
square, and has $8$ elements.
Apparently $*$, $\downarrow$, $\Downarrow$ induce an action of $D_4$
on the set of all Leonard systems in $\cal A$.
Two Leonard systems will be called {\em relatives} whenever they are
in the same orbit of this $D_4$ action. 
The relatives of $\Phi$ are as follows:

\medskip
\noindent
\begin{center}
\begin{tabular}{c|c}
name  &  relative \\
\hline
$\Phi$ & 
       $(A; \{E_i\}_{i=0}^d; A^*;  \{E^*_i\}_{i=0}^d)$ \\ 
$\Phi^{\downarrow}$ &
       $(A; \{E_i\}_{i=0}^d; A^*;  \{E^*_{d-i}\}_{i=0}^d)$ \\ 
$\Phi^{\Downarrow}$ &
       $(A; \{E_{d-i}\}_{i=0}^d; A^*;  \{E^*_i\}_{i=0}^d)$ \\ 
$\Phi^{\downarrow \Downarrow}$ &
       $(A; \{E_{d-i}\}_{i=0}^d; A^*;  \{E^*_{d-i}\}_{i=0}^d)$ \\ 
$\Phi^{*}$  & 
       $(A^*; \{E^*_i\}_{i=0}^d; A;  \{E_i\}_{i=0}^d)$ \\ 
$\Phi^{\downarrow *}$ &
       $(A^*; \{E^*_{d-i}\}_{i=0}^d; A;  \{E_i\}_{i=0}^d)$ \\ 
$\Phi^{\Downarrow *}$ &
       $(A^*; \{E^*_i\}_{i=0}^d; A;  \{E_{d-i}\}_{i=0}^d)$ \\ 
$\Phi^{\downarrow \Downarrow *}$ &
       $(A^*; \{E^*_{d-i}\}_{i=0}^d; A;  \{E_{d-i}\}_{i=0}^d)$
\end{tabular}
\end{center}
We will use the following notational convention.

\medskip

\begin{definition}
For $g \in D_4$ and for an object $f$ associated with $\Phi$ we let
$f^g$ denote the corresponding object associated with $\Phi^{g^{-1}}$.
\end{definition}

\section{The antiautomorphism $\dagger$}

Associated with the Leonard system $\Phi$, there is a certain
antiautomorphism of $\cal A$ denoted by $\dagger$ and defined as follows.
Recall that an {\em antiautomorphism} of $\cal A$
is a map $\gamma : {\cal A} \to {\cal A}$
which is an isomorphism of $\mathbb{K}$-vector spaces
and $(XY)^\gamma=Y^\gamma X^\gamma$ for all $X,Y\in {\cal A}$.

\medskip

\begin{theorem}    \cite[Theorem 7.1, Lemma 7.3]{T:qRacah}       \label{thm:anti}
Referring to the Leonard system $\Phi$,
there exists a unique antiautomorphism $\dagger$ of $\cal A$ such that
$A^{\dagger}=A$ and $A^{*\dagger}=A^*$.
Moreover $X^{\dagger\dagger}=X$ for all $X \in {\cal A}$ and
\begin{equation}     \label{eq:antiEi}
  E^{\dagger}_i= E_i, \qquad E^{*\dagger}_i= E^*_i
    \qquad (0 \leq i \leq d).
\end{equation}
\end{theorem}

\section{The parameter array}

In this section we recall some parameters associated with the
Leonard system $\Phi$.

\medskip

\begin{definition}
Referring to the Leonard system $\Phi$,
for $0 \leq i \leq d$ let $\theta_i$ 
denote the eigenvalue of $A$ associated with
$E_i$, and let  $\theta^*_i$ denote the eigenvalue
of $A^*$ associated with $E^*_i$.
We call $\theta_0, \theta_1, \ldots, \theta_d$ 
the {\em eigenvalue sequence} of $\Phi$, and
$\theta^*_0, \theta^*_1, \ldots, \theta^*_d$
the {\em dual eigenvalue sequence} of $\Phi$.
\end{definition}

\medskip

Let $V$ denote a simple $\cal A$-module.
By a {\em decomposition of $V$} we mean a sequence consisting of $1$-dimensional 
subspaces of $V$ whose direct sum is $V$.
For $0 \leq i \leq d$ we define
\begin{equation}         \label{eq:defUi}
  U_{i} = (E^*_0V+E^*_1V+\cdots+E^*_iV)\cap(E_iV+E_{i+1}V+\cdots+E_dV).
\end{equation}
By \cite[Theorem 20.7]{T:survey}
the sequence $U_0,U_1,\ldots,U_d$ is a decomposition of $V$.
We call this the {\em $\Phi$-split decomposition} of $V$.
By \cite[Theorem 20.7, Lemma 20.9]{T:survey} we have
\begin{equation}        \label{eq:AUi}
   (A-\theta_i I)U_i = U_{i+1} \quad (0 \leq i \leq d-1),
   \qquad (A-\theta_d I) U_d=0,
\end{equation}
\begin{equation}       \label{eq:AsUi}
   (A^*-\theta^*_i I)U_i = U_{i-1} \quad (1 \leq i \leq d),
   \qquad (A^*-\theta^*_0 I) U_0=0.
\end{equation}
Comparing (\ref{eq:AUi}), (\ref{eq:AsUi}) we find that for $1 \leq i \leq d$
the space $U_i$ is an eigenspace for $(A-\theta_{i-1}I)(A^*-\theta^*_iI)$,
and the corresponding eigenvalue is nonzero.
We denote this eigenvalue by $\varphi_i$. We call the sequence
$\varphi_1,\varphi_2,\ldots,\varphi_d$ the {\em first split sequence}
of $\Phi$.
The first split sequence of $\Phi^{\Downarrow}$
is denoted by $\phi_1,\phi_2,\ldots,\phi_d$, and called the
{\em second split sequence} of $\Phi$.
We emphasize that $\varphi_i \neq 0$, $\phi_i \neq 0$ for $1 \leq i \leq d$.

\medskip

\begin{definition}
By the {\em parameter array} of $\Phi$ we mean the sequence
\begin{equation}      \label{eq:param}
      (\{\theta_i\}_{i=0}^d; \{\theta^*_i\}_{i=0}^d;
        \{\varphi_i\}_{i=1}^d; \{\phi_i\}_{i=1}^d),
\end{equation}
where $\theta_0, \theta_1, \ldots, \theta_d$ 
(resp. $\theta^*_0, \theta^*_1, \ldots, \theta^*_d$)
denotes the eigenvalue sequence (resp. dual eigenvalue sequence) of $\Phi$, 
and $\varphi_1, \varphi_2, \ldots, \varphi_d$ 
(resp. $\phi_1, \phi_2, \ldots, \phi_d$)
denotes the first split sequence (resp. second split sequence) of $\Phi$.
\end{definition}

\medskip

The $D_4$ action affects the parameter array as follows.

\medskip

\begin{lemma}   \cite[Theorem 1.11]{T:Leonard}        \label{lem:D4}
Referring to the Leonard system $\Phi$,
the following (i)--(iii) hold.
\begin{itemize}
\item[(i)] The parameter array of $\Phi^*$ is 
\[
   (\{\theta^*_i\}_{i=0}^d; \{\theta_i\}_{i=0}^d;
        \{\varphi_i\}_{i=1}^d; \{\phi_{d-i+1}\}_{i=1}^d).
\]
\item[(ii)]
The parameter array of $\Phi^{\downarrow}$ is
\[
  (\{\theta_i\}_{i=0}^d; \{\theta^*_{d-i}\}_{i=0}^d;
        \{\phi_{d-i+1}\}_{i=1}^d; \{\varphi_{d-i+1}\}_{i=1}^d).
\]
\item[(iii)]
The parameter array of $\Phi^{\Downarrow}$ is
\[
  (\{\theta_{d-i}\}_{i=0}^d; \{\theta^*_{i}\}_{i=0}^d;
        \{\phi_{i}\}_{i=1}^d; \{\varphi_{i}\}_{i=1}^d).
\]
\end{itemize}
\end{lemma}

\section{Some formulae}

In this section we prove some formulae which we will find useful.
We use the following notation.

\medskip

\begin{definition}
Let $\lambda$ denote an indeterminate and let 
$\mathbb{K}[\lambda]$ denote the $\mathbb{K}$-algebra consisting of all 
polynomials in $\lambda$ that have coefficients in $\mathbb{K}$.
For $0 \leq i \leq d$  let $\tau_i$, $\eta_i$, $\tau^*_i$,
$\eta^*_i$ denote the following polynomials in $\mathbb{K}[\lambda]$:
\begin{eqnarray*}
 \tau_i &=& (\lambda-\theta_0)(\lambda-\theta_1) \cdots(\lambda-\theta_{i-1}), \\
 \eta_i &=& 
   (\lambda-\theta_d)(\lambda-\theta_{d-1})\cdots (\lambda-\theta_{d-i+1}), \\
 \tau^*_i &=& (\lambda-\theta^*_0)(\lambda-\theta^*_1) \cdots(\lambda-\theta^*_{i-1}),\\
 \eta^*_i &=& 
   (\lambda-\theta^*_d)(\lambda-\theta^*_{d-1})\cdots (\lambda-\theta^*_{d-i+1}).
\end{eqnarray*}
We observe that each of $\tau_i$, $\eta_i$, $\tau^*_i$, $\eta^*_i$ is monic
of degree $i$.
\end{definition}

\begin{theorem}           \label{thm:rel}
Referring to the Leonard system $\Phi$,
for $0 \leq i \leq d$ we have
\begin{eqnarray}
\eta_i(A)E^*_0E_0 &=&
 \frac{\phi_{1}\phi_{2} \cdots \phi_{i}}
      {\eta^*_d(\theta^*_0)} 
         \eta^*_{d-i}(A^*)E_0, 
      \label{eq:rel1}  \\
\eta_i(A)E^*_dE_0 &=&
 \frac{\varphi_d \varphi_{d-1} \cdots \varphi_{d-i+1}}
      {\tau^*_d (\theta^*_d)}
           \tau^*_{d-i}(A^*)E_0, 
       \label{eq:reld}  \\
\tau_i(A) E^*_0 E_d &=&
 \frac{\varphi_{1}\varphi_{2} \cdots \varphi_{i}}
      {\eta^*_d(\theta^*_0)} 
         \eta^*_{d-i}(A^*)E_d, 
           \label{eq:relD}   \\
\tau_i(A) E^*_d E_d &=&
 \frac{\phi_{d}\phi_{d-1} \cdots \phi_{d-i+1}}
      {\tau^*_d(\theta^*_d)} 
         \tau^*_{d-i}(A^*)E_d,  
             \label{eq:reldD}\\
\eta^*_i(A^*) E_0 E^*_0 &=&
 \frac{\phi_{d}\phi_{d-1} \cdots \phi_{d-i+1}}
      {\eta_d(\theta_0)} 
         \eta_{d-i}(A)E^*_0,
         \label{eq:rels} \\
\eta^*_i(A^*) E_d E^*_0 &=&
 \frac{\varphi_{d}\varphi_{d-1} \cdots \varphi_{d-i+1}}
      {\tau_d(\theta_d)} 
         \tau_{d-i}(A)E^*_0, 
           \label{eq:relds} \\
\tau^*_i(A^*) E_0 E^*_d &=&
 \frac{\varphi_{1}\varphi_{2} \cdots \varphi_{i}}
      {\eta_d(\theta_0)} 
         \eta_{d-i}(A)E^*_d, 
            \label{eq:relDs} \\
\tau^*_i(A^*) E_d E^*_d &=&
 \frac{\phi_{1}\phi_{2} \cdots \phi_{i}}
      {\tau_d(\theta_d)} 
         \tau_{d-i}(A)E^*_d.
             \label{eq:reldDs}
\end{eqnarray}
\end{theorem}

\begin{proof}
We first prove (\ref{eq:relD}).
Let $i$ be given and define
\[
   H = \eta^*_d(\theta^*_0) \tau_i(A)E^*_0
            - \varphi_1 \varphi_2 \cdots \varphi_i \eta^*_{d-i}(A^*).
\]
We show $H E_d=0$. To this end we let $V$ denote a simple
$\cal A$-module and show $H E_d V=0$.
By (\ref{eq:defEi}) we have
\[
      \eta^*_d(\theta^*_0) E^*_0 = \eta^*_d(A^*).
\]
Using this we find that $H$ is equal to
\begin{equation}         \label{eq:relaux1}
 (A-\theta_{i-1}I)\cdots(A-\theta_0I)(A^*-\theta^*_1I)\cdots(A^*-\theta^*_iI)
       - \varphi_1 \cdots \varphi_i I
\end{equation}
times
\[
     \eta^*_{d-i}(A^*).
\]
Setting $i=d$ in (\ref{eq:defUi}) we find $E_dV=U_d$.
Recursively using (\ref{eq:AsUi}) we find
\[
   \eta^*_{d-i}(A^*)U_d=U_i.
\]
By (\ref{eq:AUi}), (\ref{eq:AsUi}) and the definition of $\varphi_i$,
we find that (\ref{eq:relaux1})
vanishes on $U_i$.
By these comments we find $H E_d V=0$ as desired. Thus (\ref{eq:relD})
holds.
To obtain the remaining equations apply $D_4$
and use Lemma \ref{lem:D4}.
\end{proof}

\medskip

Applying the antiautomorphism $\dagger$ to (\ref{eq:rel1})--(\ref{eq:reldDs})
we routinely obtain the following result.

\medskip

\begin{corollary}           \label{cor:relanti}
Referring to the Leonard system $\Phi$,
for $0 \leq i \leq d$ we have
\begin{eqnarray}
E_0E^*_0\eta_i(A) &=&
 \frac{\phi_{1}\phi_{2} \cdots \phi_{i}}
      {\eta^*_d(\theta^*_0)} 
         E_0 \eta^*_{d-i}(A^*), 
      \label{eq:rel1anti}  \\
E_0E^*_d \eta_i(A) &=&
 \frac{\varphi_d \varphi_{d-1} \cdots \varphi_{d-i+1}}
      {\tau^*_d (\theta^*_d)}
           E_0 \tau^*_{d-i}(A^*), 
       \label{eq:reldanti}  \\
E_d E^*_0 \tau_i(A) &=&
 \frac{\varphi_{1}\varphi_{2} \cdots \varphi_{i}}
      {\eta^*_d(\theta^*_0)} 
         E_d \eta^*_{d-i}(A^*), 
           \label{eq:relDanti}   \\
E_d E^*_d \tau_i(A) &=&
 \frac{\phi_{d}\phi_{d-1} \cdots \phi_{d-i+1}}
      {\tau^*_d(\theta^*_d)} 
         E_d \tau^*_{d-i}(A^*),  
             \label{eq:reldDanti}\\
E^*_0 E_0 \eta^*_i(A^*) &=&
 \frac{\phi_{d}\phi_{d-1} \cdots \phi_{d-i+1}}
      {\eta_d(\theta_0)} 
         E^*_0 \eta_{d-i}(A),
         \label{eq:relsanti} \\
E^*_0 E_d \eta^*_i(A^*)  &=&
 \frac{\varphi_{d}\varphi_{d-1} \cdots \varphi_{d-i+1}}
      {\tau_d(\theta_d)} 
         E^*_0 \tau_{d-i}(A), 
           \label{eq:reldsanti} \\
E^*_d E_0 \tau^*_i(A^*) &=&
 \frac{\varphi_{1}\varphi_{2} \cdots \varphi_{i}}
      {\eta_d(\theta_0)} 
         E^*_d \eta_{d-i}(A), 
            \label{eq:relDsanti} \\
E^*_d E_d \tau^*_i(A^*)  &=&
 \frac{\phi_{1}\phi_{2} \cdots \phi_{i}}
      {\tau_d(\theta_d)} 
         E^*_d \tau_{d-i}(A).
             \label{eq:reldDsanti}
\end{eqnarray}
\end{corollary}

\medskip

To obtain further formulae concerning $\Phi$
we need the following general lemma.

\medskip

\begin{lemma}      \label{lem:aux}
Let $\xi_0,\xi_1,\ldots,\xi_d$ denote a sequence
of mutually distinct scalars in $\mathbb{K}$.
Then for an indeterminate $\lambda$,
\begin{equation}       \label{eq:aux}
 \sum_{i=0}^{d} 
  \frac{(\lambda-\xi_0)(\lambda-\xi_1)\cdots(\lambda-\xi_{i-1})}
       {(\xi_0-\xi_1)(\xi_0-\xi_2)\cdots(\xi_0-\xi_i)}
= \frac{(\lambda-\xi_1)(\lambda-\xi_2)\cdots(\lambda-\xi_{d})}
       {(\xi_0-\xi_1)(\xi_0-\xi_2)\cdots(\xi_0-\xi_d)}.
\end{equation}
\end{lemma}

\begin{proof}
We will use induction on $d$. Assume $d \geq 1$; otherwise the result
is trivial.
The left side of (\ref{eq:aux}) is equal to
\begin{eqnarray*}
 & & 
 \sum_{i=0}^{d-1} 
  \frac{(\lambda-\xi_0)(\lambda-\xi_1)\cdots(\lambda-\xi_{i-1})}
       {(\xi_0-\xi_1)(\xi_0-\xi_2)\cdots(\xi_0-\xi_i)}
 + \frac{(\lambda-\xi_0)(\lambda-\xi_1)\cdots(\lambda-\xi_{d-1})}
       {(\xi_0-\xi_1)(\xi_0-\xi_2)\cdots(\xi_0-\xi_d)} \\
&=& \frac{(\lambda-\xi_1)(\lambda-\xi_2)\cdots(\lambda-\xi_{d-1})}
         {(\xi_0-\xi_1)(\xi_0-\xi_2)\cdots(\xi_0-\xi_{d-1})} 
 +\frac{(\lambda-\xi_0)(\lambda-\xi_1)\cdots(\lambda-\xi_{d-1})}
       {(\xi_0-\xi_1)(\xi_0-\xi_2)\cdots(\xi_0-\xi_d)} 
    \qquad (\text{by induction}) \\
&=& \frac{(\lambda-\xi_1)(\lambda-\xi_2)\cdots(\lambda-\xi_{d-1})}
         {(\xi_0-\xi_1)(\xi_0-\xi_2)\cdots(\xi_0-\xi_{d-1})} 
    \frac{\lambda-\xi_d}
         {\xi_0-\xi_d}
\end{eqnarray*}
and this is equal to the right side of (\ref{eq:aux}).
\end{proof}

\medskip

\begin{proposition}     \label{thm:trasision}
The following hold.
\begin{equation}
\tau_d = \sum_{i=0}^{d} \tau_{d-i}(\theta_d) \eta_i, \qquad \qquad
\eta_d = \sum_{i=0}^{d} \eta_{d-i}(\theta_0) \tau_i, \label{eq:transetad}
\end{equation}
\begin{equation}
\tau^*_d = \sum_{i=0}^{d} \tau^*_{d-i}(\theta^*_d) \eta^*_i, \qquad \qquad
\eta^*_d = \sum_{i=0}^{d} \eta^*_{d-i}(\theta^*_0) \tau^*_i. 
\label{eq:transetasd}
\end{equation}
\end{proposition}

\begin{proof}
In  (\ref{eq:aux}) we set $\xi_i=\theta_i$ for $0 \leq i \leq d$
and multiply both sides by $\eta_d(\theta_0)$;
this gives the equation on the right in (\ref{eq:transetad}).
To obtain the remaining equations apply $D_4$
and use Lemma \ref{lem:D4}.
\end{proof}

\medskip

\begin{theorem}      \label{thm:E0EsdEdEs0}
Referring to the Leonard system $\Phi$,
the following hold.
\begin{eqnarray}
E_0 E^*_d E_d E^*_0 &=&
  \frac{\varphi_1 \varphi_2 \cdots \varphi_d}
       {\tau_d(\theta_d) \tau^*_d(\theta^*_d)} E_0 E^*_0,  
     \label{eq:E0EsdEdEs0} \\
E_0 E^*_0 E_d E^*_d &=&
  \frac{\phi_1 \phi_2 \cdots \phi_d}
       {\tau_d(\theta_d) \eta^*_d(\theta^*_0)} E_0 E^*_d,
     \label{eq:E0Es0EdEsd}  \\
E_d E^*_d E_0 E^*_0 &=&
  \frac{\phi_1 \phi_2 \cdots \phi_d}
       {\eta_d(\theta_0) \tau^*_d(\theta^*_d)} E_d E^*_0,
     \label{eq:EdEsdE0Es0}  \\
E_d E^*_0 E_0 E^*_d &=&
  \frac{\varphi_1 \varphi_2 \cdots \varphi_d}
       {\eta_d(\theta_0) \eta^*_d(\theta^*_0)} E_d E^*_d,
     \label{eq:EdEs0E0Esd}  \\
E^*_0 E_d E^*_d E_0 &=&
  \frac{\varphi_1 \varphi_2 \cdots \varphi_d}
       {\tau_d(\theta_d) \tau^*_d(\theta^*_d)} E^*_0 E_0,  
     \label{eq:Es0EdEsdE0}  \\
E^*_0 E_0 E^*_d E_d &=&
  \frac{\phi_1 \phi_2 \cdots \phi_d}
       {\eta_d(\theta_0) \tau^*_d(\theta^*_d)} E^*_0 E_d,
     \label{eq:Es0E0EsdEd}  \\
E^*_d E_d E^*_0 E_0 &=&
  \frac{\phi_1 \phi_2 \cdots \phi_d}
       {\tau_d(\theta_d) \eta^*_d(\theta^*_0)} E^*_d E_0,
     \label{eq:EsdEdEs0E0}  \\
E^*_d E_0 E^*_0 E_d &=&
  \frac{\varphi_1 \varphi_2 \cdots \varphi_d}
       {\eta_d(\theta_0) \eta^*_d(\theta^*_0)} E^*_d E_d.
     \label{eq:EsdE0Es0Ed}
\end{eqnarray}
\end{theorem}

\begin{proof}
We first prove (\ref{eq:E0EsdEdEs0}).
From (\ref{eq:defEi}) we have
\[
   E_d = \frac{1}{\tau_d(\theta_d)} \tau_d(A).
\]
In this equation we evaluate $\tau_d(A)$ using the
equation on the left in (\ref{eq:transetad}) to get
\[
  E_d = \frac{1}{\tau_d(\theta_d)} \sum_{i=0}^{d} \tau_{d-i}(\theta_d)\eta_i(A).
\]
In this equation we multiply each term on the left by
$E_0E^*_d$ and on the right by $E^*_0$; the result is
\[
  E_0 E^*_d E_d E^*_0
 = \frac{1}{\tau_d(\theta_d)}
     \sum_{i=0}^{d} \tau_{d-i}(\theta_d) E_0 E^*_d \eta_i(A) E^*_0.
\]
In this equation we evaluate each term $E_0 E^*_d \eta_i(A)$ 
using (\ref{eq:reldanti}); we find $E_0E^*_dE_dE^*_0$ is equal to
\[
    \frac{\varphi_1\varphi_2 \cdots\varphi_d}
         {\tau_d(\theta_d)\tau^*_d(\theta^*_d)}
\]
times
\begin{equation}            \label{eq:E0EsdEdEs0aux}
    \sum_{i=0}^{d} 
      \frac{\tau_{d-i}(\theta_d)  E_0 \tau^*_{d-i}(A^*)E^*_0}
           {\varphi_1 \varphi_2 \cdots \varphi_{d-i}}.
\end{equation}
We show (\ref{eq:E0EsdEdEs0aux}) is equal to $E_0E^*_0$.
Recall $A^*E^*_0 = \theta^*_0E^*_0$, so
$\tau^*_j(A^*)E^*_0 = \tau^*_j(\theta^*_0)E^*_0$ for $0 \leq j \leq d$.
Also $\tau^*_j(\theta^*_0)$ is equal to $1$ for $j=0$ and $0$
for $1 \leq j \leq d$.
Combining these comments we find
$\tau^*_j(A^*)E^*_0= \delta_{0,j} E^*_0$ for $0 \leq j \leq d$.
Evaluating (\ref{eq:E0EsdEdEs0aux}) using this we find it is equal to
$E_0E^*_0$. Thus (\ref{eq:E0EsdEdEs0}) holds.
To obtain the remaining equations apply $D_4$
and use Lemma \ref{lem:D4}.
\end{proof}

\section{The $\Phi$-split basis}

Let $V$ denote a simple $\cal A$-module.
In this section we associate with $\Phi$ a certain basis for $V$,
and we consider the action of $\cal A$ on this basis.

Consider the decomposition $U_0,U_1,\ldots,U_d$ from (\ref{eq:defUi}).
For $0 \leq i \leq d$ let $u_i$ denote a nonzero vector in $U_i$,
and observe $u_0,u_1,\ldots,u_d$ is a basis for $V$.
In view of (\ref{eq:AUi}) we normalize this basis so that
\begin{equation}      \label{eq:Aui}
   (A-\theta_iI) u_i=u_{i+1} \quad (0 \leq i \leq d-1),   \qquad
   (A-\theta_d I)u_d = 0.
\end{equation}
By (\ref{eq:AsUi}), (\ref{eq:Aui}) and the definition of $\varphi_i$,
we find
\begin{equation}     \label{eq:Asui}
   (A^* - \theta^*_iI) u_i = \varphi_i u_{i-1} \quad (1 \leq i \leq d),  
   \qquad
   (A^* - \theta^*_0I) u_0 = 0.
\end{equation}
We call $u_0,u_1,\ldots,u_d$ the {\em $\Phi$-split basis} for $V$.
The $\Phi$-split basis is unique up to multiplication
of each basis vector by the same nonzero scalar in $\mathbb{K}$.

\medskip

\begin{definition}          \label{def:natural}
For $X \in {\cal A}$, let $X^{\natural}$ denote the matrix in $\Mat{d+1}$
that represents $X$ with respect to the $\Phi$-split basis.
We observe $\natural: {\cal A} \to \Mat{d+1}$ is an isomorphism
of $\mathbb{K}$-algebras.
\end{definition}

\begin{example}        \label{exm:AAs}
We have
\[
A^{\natural} =
\begin{pmatrix}
   \theta_0 &          & & & & \text{\bf 0} \\
   1        & \theta_1 \\
            & 1        & \theta_2 \\
            &          &  \cdot   &  \cdot \\
            &          &          & \cdot  & \cdot\\
   \text{\bf 0} &      &          &       & 1 & \theta_d
\end{pmatrix},
\qquad
A^{*\natural} =
\begin{pmatrix}
   \theta^*_0 &  \varphi_1  & & & & \text{\bf 0} \\
           & \theta^*_1 & \varphi_2 \\
            &         & \theta^*_2 & \cdot \\
            &          &       &  \cdot & \cdot \\
            &          &          &   & \cdot & \varphi_d\\
   \text{\bf 0} &      &          &       &  & \theta^*_d
\end{pmatrix}.
\]
\end{example}

\begin{proof}
Immediate from  (\ref{eq:Aui}), (\ref{eq:Asui}).
\end{proof}

\begin{definition}      \label{def:Delta}
Referring to Definition \ref{def:natural},
for $0 \leq i,j \leq d$ let $\Delta_{i,j}$ denote the element of 
$\cal A$ such that $(\Delta_{i,j})^{\natural}$ has $(i,j)$-entry $1$ and
all other entries $0$. 
\end{definition}

\medskip

We have a comment.

\medskip

\begin{lemma}          \label{lem:matrixunits}
Referring to Definition \ref{def:Delta},
the $\mathbb{K}$-vector space $\cal A$ has basis
\[
   \Delta_{i,j} \qquad\qquad  (0 \leq i,j \leq d).
\]
Moreover
\begin{equation}
   \Delta_{i,j} \Delta_{r,s} = \delta_{j,r} \Delta_{i,s}
    \qquad (0 \leq i,j,r,s \leq d).
\end{equation}
\end{lemma}

\medskip

We are going to give some formulae for the $\Delta_{i,j}$.
We will use the following notation.
Recall by  \cite[Lemma 10.2]{T:qRacah} that the trace of $E^*_0E_0$
is nonzero; let $\nu$ denote the multiplicative inverse of this trace,
so that
\[
  \text{\rm tr}(E^*_0E_0)=\nu^{-1}.
\]
By \cite[Theorem 18.8]{T:qRacah} we have
\begin{equation}    \label{eq:nu}
\nu =  
   \frac{\eta_d(\theta_0)\eta^*_d(\theta^*_0)}
        {\phi_1\phi_2\cdots\phi_d}.
\end{equation}
Applying $D_4$ to (\ref{eq:nu}) and using Lemma \ref{lem:D4} we find
\begin{equation}        \label{eq:nu2}
\nu^{\downarrow\Downarrow} =
   \frac{\tau_d(\theta_d)\tau^*_d(\theta^*_d)}
        {\phi_1\phi_2\cdots\phi_d}.
\end{equation}

\medskip

\begin{theorem}         \label{thm:main}
Referring to Definition \ref{def:Delta},
$\Delta_{i,j}$ is equal to each of
\begin{equation}        \label{eq:Deltaij}
     \frac{\nu \tau_i(A) E^*_0 E_0 \tau^*_j(A^*)}
          {\varphi_1 \varphi_2 \cdots \varphi_j},
  \qquad \qquad
     \frac{\nu^{\downarrow\Downarrow} \eta^*_{d-i}(A^*) E_d E^*_d \eta_{d-j}(A)}
          {\varphi_{d} \varphi_{d-1} \cdots \varphi_{i+1}}
\end{equation}
for $0 \leq i,j \leq d$.
\end{theorem}

\begin{proof}
We first prove that $\Delta_{i,j}$
is equal to the expression on the left in (\ref{eq:Deltaij}).
Let $V$ denote a simple $\cal A$-module and let
$u_0,u_1,\ldots,u_d$ denote a $\Phi$-split basis for $V$.
By Definition \ref{def:Delta} it suffices to show
\begin{equation}      \label{eq:mainaux1}
  \nu \tau_i(A) E^*_0 E_0 \tau^*_j(A^*) u_k =
    \delta_{j,k} \varphi_1 \varphi_2 \cdots \varphi_j u_i
             \qquad (0 \leq k  \leq d).
\end{equation}
By a routine induction using (\ref{eq:Asui}) we find that, 
when we write $\tau^*_j(A^*) u_k$ as a linear combination of 
$u_0,u_1,\ldots,u_d$, the coefficient of
$u_0$ is $\delta_{j,k}\varphi_1\varphi_2 \cdots \varphi_j$.
We claim that
\begin{equation}    \label{eq:mainaux2}
E^*_0 E_0 u_t = \delta_{0,t} \nu^{-1} u_0   \qquad (0 \leq t \leq d).
\end{equation}
To prove (\ref{eq:mainaux2}) we consider the matrix 
$(E^*_0E_0)^\natural$.
For $1 \leq r \leq d$, column $r$ of 
 $(E^*_0E_0)^\natural$
is zero since $E_0u_r = 0$ by (\ref{eq:defUi}).
For $1 \leq s \leq d$, row $s$ of 
$(E^*_0E_0)^\natural$
is zero since $E^*_0V = U_0$ and since $u_0$ is a basis
for $U_0$.
Apparently the $(0,0)$-entry of 
$(E^*_0E_0)^\natural$
is equal to the trace of $E^*_0E_0$ or in other words $\nu^{-1}$.
Line (\ref{eq:mainaux2}) follows and the claim is proved.
From (\ref{eq:Aui}) we find $\tau_i(A)u_0 = u_i$.
By these commments we routinely obtain (\ref{eq:mainaux1}). 
We have now shown that $\Delta_{i,j}$ is equal to the expression
on the left in (\ref{eq:Deltaij}). 

Next we prove that $\Delta_{i,j}$ is equal to the expression on the right
in (\ref{eq:Deltaij}).
In the expression on the left in (\ref{eq:Deltaij})
 we evaluate $\tau_i(A)E^*_0$ using (\ref{eq:relds})
and evaluate $E_0 \tau^*_j(A^*)$ using (\ref{eq:reldanti});
the result is
\[
\Delta_{i,j} =
  \frac{\nu \tau_d(\theta_d)\tau^*_d(\theta^*_d)}
       {\varphi_1\varphi_2\cdots\varphi_d \cdot 
        \varphi_d \varphi_{d-1}\cdots\varphi_{i+1}}
    \eta^*_{d-i}(A^*)E_d E^*_0 E_0 E^*_d \eta_{d-j}(A).
\]
In this equation we evaluate $E_d E^*_0 E_0 E^*_d$ using 
(\ref{eq:EdEs0E0Esd}); the result is
\begin{equation}           \label{eq:aux1}
\Delta_{i,j} =
  \frac{\nu \tau_d(\theta_d)\tau^*_d(\theta^*_d)}
       {\eta_d(\theta_0)\eta^*_d(\theta^*_0)
             \varphi_d \varphi_{d-1}\cdots\varphi_{i+1}}
    \eta^*_{d-i}(A^*)E_d E^*_d \eta_{d-j}(A).
\end{equation}
From (\ref{eq:nu}) and (\ref{eq:nu2}),
\begin{equation}      \label{eq:aux2}
    \frac{\nu \tau_d(\theta_d)\tau^*_d(\theta^*_d)}
       {\eta_d(\theta_0)\eta^*_d(\theta^*_0)} = \nu^{\downarrow\Downarrow}.
\end{equation}
Combining (\ref{eq:aux1}), (\ref{eq:aux2}) we find that $\Delta_{i,j}$
is equal to the expression on the right in (\ref{eq:Deltaij}).
\end{proof}

\section{The primitive idempotents associated with the $\Phi$-split decomposition}

In this section we consider the following elements in $\cal A$.

\medskip

\begin{definition}      \label{def:Fi}
Referring to Definition \ref{def:Delta}, 
we define
\[
      F_i=\Delta_{i,i}  \qquad (0 \leq i \leq d).
\]
\end{definition}

\medskip

We have two comments.

\begin{lemma}
Let $V$ denote a simple $\cal A$-module.
Then for $0 \leq i \leq d$, $F_i$ satisfies
\[
    (F_i-I)U_i=0,
\]
\[
  F_i U_j = 0 \quad \text{ if $i \neq j$ $(0 \leq j \leq d)$}.
\]
In other words $F_i$ is the projection from $V$ onto $U_i$.
\end{lemma}

\begin{lemma}
The following (i)--(iii) hold.
\begin{itemize}
\item[(i)] $F_iF_j = \delta_{i,j} F_i$ $(0 \leq i,j \leq d)$.
\item[(ii)] $F_0+F_1+\cdots+F_d = I$.
\item[(iii)] Each of $F_0,F_1,\ldots,F_d$ has rank $1$.
\end{itemize}
\end{lemma}

\medskip

\begin{corollary}           \label{cor:Frelative}
For $0 \leq i \leq d$, $F_i$ is equal to each of
\begin{equation}        \label{eq:Fi}
     \frac{\nu \tau_i(A) E^*_0 E_0 \tau^*_i(A^*)}
          {\varphi_1 \varphi_2 \cdots \varphi_i},
  \qquad \qquad
     \frac{\nu^{\downarrow\Downarrow} \eta^*_{d-i}(A^*) E_d E^*_d \eta_{d-i}(A)}
          {\varphi_{d} \varphi_{d-1} \cdots \varphi_{i+1}}.
\end{equation}
\end{corollary}

\begin{proof}
Combine Definition \ref{def:Fi} and Theorem \ref{thm:main}.
\end{proof}

\section{Some formulae for the split sequences}

In this section we give some formulae for the first and second
split sequences of $\Phi$ in terms of $F_i$, $F^*_i$,
$F^{\downarrow}_i$, $F^{\Downarrow}_i$.
We start with a lemma.

\medskip

\begin{lemma}     \label{lem:comFk}
For $0 \leq k \leq d$,
\begin{equation}      \label{eq:comFk}
  {\rm tr}((AA^* - A^*A)F_k) = \varphi_k - \varphi_{k+1},
\end{equation}
where $\varphi_0=0$ and $\varphi_{d+1}=0$.
\end{lemma}

\begin{proof}
The left side of (\ref{eq:comFk}) is equal to the trace of 
$(A^{\natural}A^{*\natural}- A^{*\natural}A^{\natural})F^{\natural}_k$.
By Definition \ref{def:Fi}, $F^{\natural}_k$ has
$(k,k)$-entry $1$ and all other entries $0$.
Using this we find that for $B \in \Mat{d+1}$, 
the trace of $B F^{\natural}_k$ is equal to the $(k,k)$-entry of $B$.
By Example \ref{exm:AAs} the $(k,k)$-entry of
$A^{\natural} A^{*\natural} - A^{*\natural}A^{\natural}$ is equal to
$\varphi_{k}-\varphi_{k+1}$.
By these comments we obtain (\ref{eq:comFk}).
\end{proof}

\medskip

We now come to our final theorem.

\medskip

\begin{theorem}          \label{thm:traceformula}
For $1 \leq i \leq d$, $\varphi_i$ is equal to each of
\[
 \sum_{k=i}^{d} {\rm tr} ((AA^*-A^*A)F_k), 
\qquad 
 \sum_{k=0}^{i-1} {\rm tr} ((A^*A-AA^*)F_k),
\]
\[
  \sum_{k=0}^{i-1} {\rm tr} ((AA^*-A^*A)F^*_k),
 \qquad
  \sum_{k=i}^{d} {\rm tr} ((A^*A-AA^*)F^*_k).
\]
Moreover $\phi_i$ is equal to each of
\[
  \sum_{k=d-i+1}^{d} {\rm tr} ((AA^*-A^*A)F^{\downarrow}_k),
 \qquad 
  \sum_{k=0}^{d-i} {\rm tr} ((A^*A-AA^*)F^{\downarrow}_k),
\] 
\[
  \sum_{k=i}^{d} {\rm tr} ((AA^*-A^*A)F^{\Downarrow}_k),
 \qquad
  \sum_{k=0}^{i-1} {\rm tr} ((A^*A-AA^*)F^{\Downarrow}_k).
\]
\end{theorem}

\begin{proof}
Referring to the first displayed line of the theorem,
we find $\varphi_i$ is equal to the expression
on the left (resp. right) by summing (\ref{eq:comFk}) over $k=i,\ldots,d$
(resp. $k=0,\ldots,i-1$).
To obtain the remaining equations apply $D_4$
and use Lemma \ref{lem:D4}.
\end{proof}

\bigskip

\bibliographystyle{plain}

\begin{thebibliography}{1}

\bibitem{AC}
H.~Alnajjar  and B.~Curtin.
\newblock
A family of tridiagonal pairs.
\newblock {\em
Linear Algebra Appl.}
{\bf 390}
(2004)
369--384.

\bibitem{AC2}
H.~Alnajjar  and B.~Curtin.
\newblock
A family of tridiagonal pairs related to
the quantum affine algebra 
$U\sb q(\widehat{\mathfrak{sl}}\sb 2)$.
\newblock {\em
Electron. J. Linear Algebra}
{\bf 13} (2005) 1--9. 


\bibitem{H}
B.~Hartwig. 
Three mutually adjacent Leonard pairs. 
{\em Linear Algebra Appl.} {\bf 408} (2005) 19--39;
{\tt arXiv:math.AC/0508415}.



\bibitem{ITT}
T.~Ito, K.~Tanabe, and P.~Terwilliger.
\newblock Some algebra related to ${P}$- and ${Q}$-polynomial association
  schemes,  in:
\newblock {\em Codes and Association Schemes (Piscataway NJ, 1999)}, Amer.
Math. Soc., Providence RI, 2001, pp.
     167--192; 
{\tt arXiv:math.CO/0406556}.

\bibitem{IT:shape}
T.~Ito and P.~Terwilliger.
\newblock The shape of a tridiagonal pair.
\newblock {\em J. Pure Appl. Algebra} {\bf 188} (2004) 145--160;
{\tt arXiv:math.QA/0304244}.

\bibitem{IT:uqsl2hat}
T.~Ito and P.~Terwilliger.
\newblock {Tridiagonal pairs and the quantum affine 
algebra
$U_q({\widehat{sl}}_2)$.}
\newblock {\em Ramanujan J.} In press; 
{\tt arXiv:math.QA/0310042}.

\bibitem{IT:non-nilpotent}
T.~Ito and P.~Terwilliger.
\newblock
Two non-nilpotent linear transformations that satisfy the cubic $q$-Serre relations.
{\em J. Algebra}. Submitted;
{\tt arXiv:math.QA/0508398}.

\bibitem{ITW:equitable}
T.~Ito, P.~Terwilliger and C.~Weng.
\newblock {The quantum algebra $U_q(sl_2)$ and its equitable presentation.}
{\em J. Algebra}. In press;
{\tt arXiv:math.QA/0507477}.

\bibitem{N:aw}
K.~Nomura.
\newblock
Tridiagonal pairs and the {A}skey-{W}ilson relations.
\newblock {\em Linear Algebra Appl.}
{\bf 397} (2005) 99--106.

\bibitem{N:refine}
K.~Nomura.
\newblock A refinement of the split decomposition of
a tridiagonal pair.
\newblock {\em Linear Algebra Appl.}
{\bf 403} (2005) 1--23.

\bibitem{N:height1}
K.~Nomura.
\newblock
Tridiagonal pairs of height one.
\newblock {\em Linear Algebra Appl.}
{\bf 403} (2005) 118--142.

\bibitem{NT:balanced}
K.~Nomura and P.~Terwilliger.
\newblock
Balanced Leonard pairs.
\newblock {\em Linear Algebra Appl.}
Submitted;
{\tt arXiv:math.RA/0506219}. 

\bibitem{NT:formula}
K.~Nomura and P.~Terwilliger.
\newblock
Some trace formulae involving the split sequences of a Leonard pair.
{\em Linear Algebra Appl.}
{\bf 413} (2006), 189--201;
{\tt arXiv:math.RA/0508407}.

\bibitem{NT:det}
K.~Nomura and P.~Terwilliger.
\newblock
The determinant of $AA^*-A^*A$ for a Leonard pair $A,A^*$.
{\em Linear Algebra Appl.}
In press;
{\tt arXiv:math.RA/0511641}.

\bibitem{P}
A.~A.~Pascasio.
\newblock     
On the multiplicities of the primitive idempotents of a
 {$Q$}-polynomial distance-regular graph.
\newblock{ \em
European J. Combin.}
{\bf 23}
(2002)
1073--1078.

\bibitem{T:sub1}
P.~Terwilliger.
\newblock The subconstituent algebra of an association scheme I. 
\newblock {\em J. Algebraic Combin.} {\bf 1} (1992) 363--388.
   
\bibitem{T:sub3}
P.~Terwilliger.
\newblock The subconstituent algebra of an association scheme III.
\newblock{\em
J. Algebraic Combin. }
{\bf 2}  (1993) 177--210.

\bibitem{T:Leonard}
P.~Terwilliger.
\newblock Two linear transformations each tridiagonal with respect to an
  eigenbasis of the other.
  \newblock {\em Linear Algebra Appl.}  {\bf 330} (2001) 149--203;
{\tt arXiv:math.RA/0406555}.

\bibitem{T:qSerre}
  P.~Terwilliger.
  \newblock Two relations that generalize the $q$-Serre relations and the
  Dolan-Grady relations. In
  \newblock {\em  Physics and
  Combinatorics 1999 (Nagoya)}, 377--398, World Scientific Publishing,
   River Edge, NJ, 2001; 
{\tt arXiv:math.QA/0307016}.

\bibitem{T:24points}
   P.~Terwilliger.
   \newblock  Leonard pairs from 24 points of view.
   \newblock {\em Rocky Mountain J. Math.} {\bf 32}(2) (2002) 827--888;
{\tt arXiv:math.RA/0406577}.

\bibitem{T:conform}
   P.~Terwilliger.
   \newblock Two linear transformations each tridiagonal with respect to an
     eigenbasis of the other; the $TD$-$D$ and the $LB$-$UB$ canonical form.
\newblock {\em J. Algebra}. {\bf 291} (2005) 1--45;
{\tt arXiv:math.RA/0304077}.

\bibitem{T:intro}
    P.\ Terwilliger.
    \newblock Introduction to Leonard pairs.
    \newblock {OPSFA Rome 2001}.
    \newblock{\em J. Comput. Appl. Math.} {\bf 153}(2) (2003)
    463--475.

\bibitem{T:intro2}
P.\ Terwilliger.
Introduction to {L}eonard pairs and
  {L}eonard systems.
 {\em S\=uri\-kaiseki\-kenky\=usho} {\em K\=oky\=uroku},
 (1109) 67--79, 1999.   Algebraic combinatorics  (Kyoto, 1999).

\bibitem{T:split}
P.~Terwilliger.
\newblock Two linear transformations each tridiagonal with respect to an
  eigenbasis of the other; comments on the split decomposition.
\newblock {\em 
 J. Comput. Appl. Math.} {\bf 178} (2005) 437--452;
{\tt arXiv:math.RA/0306290}.

 \bibitem{T:array}
 P.~Terwilliger.
 \newblock Two linear transformations each tridiagonal with respect to an
   eigenbasis of the other; comments on the parameter array.
\newblock {\em
Des. Codes Cryptogr.}  {\bf 34}  (2005) 307--332;
{\tt arXiv:math.RA/0306291}.

\bibitem{T:qRacah}
P.~Terwilliger.
\newblock Leonard pairs and the $q$-Racah polynomials.
\newblock {\em Linear Algebra Appl.} {\bf 387} (2004) 235--276;
{\tt arXiv:math.QA/0306301}.

\bibitem{T:survey}
P.~Terwilliger.
\newblock
Two linear transformations each tridiagonal with respect to an eigenbasis
of the other; an algebraic approach to the Askey scheme of
orthogonal polynomials.
\newblock Lecture notes for the summer school on orthogonal polynomials
and special functions. Universidad Carlos III de Madrid, Leganes, Spain.
July 8--July 18, 2004;
{\tt arXiv:math.QA/0408390}. 

\bibitem{T:Kac-Moody}
P.~Terwilliger.
\newblock
The equitable presentation for the quantum group $U_q(g)$ 
associated with a symmetrizable Kac-Moody algebra $g$.
{\em J. Algebra}. In press; 
{\tt arXiv:math.QA/0507478}.

\bibitem{TV}
P.~Terwilliger and R.~Vidunas.
\newblock Leonard pairs and the Askey-Wilson relations.
\newblock {\em J. Algebra Appl.} {\bf 3} (2004) 411--426;
{\tt arXiv:math.QA/0305356}.

\bibitem{V}
R.~Vidunas.
\newblock
Normalized Leonard pairs and Askey-Wilson relations. Preprint;
\hfil\break
{\tt 
arXiv:math.RA/0505041}.

\end{thebibliography}

\bigskip\bigskip\noindent
Kazumasa Nomura\\
College of Liberal Arts and Sciences\\
Tokyo Medical and Dental University\\
Kohnodai, Ichikawa, 272-0827 Japan\\
email: nomura.las@tmd.ac.jp

\bigskip\noindent
Paul Terwilliger\\
Department of Mathematics\\
University of Wisconsin\\
480 Lincoln drive, Madison, Wisconsin, 53706 USA\\
email: terwilli@math.wisc.edu

\bigskip\noindent
{\bf Keywords.}
Leonard pair, tridiagonal pair, $q$-Racah polynomial, orthgonal polynomial.

\noindent
{\bf 2000 Mathematics Subject Classification}.
05E35, 05E30, 33C45, 33D45.

\end{document}